

\magnification\magstephalf
\baselineskip14pt

\def\Pfaff/{P\kern-.07emfaff}
\let\blash=\backslash
\font\sc=cmcsc10
\def\pfbox
  {\hbox{\hskip 3pt\vbox{\hrule
  \hbox to 7pt{\vrule height 7pt\hfill\vrule}
  \hrule}}\hskip3pt}
\def\bib[#1] {\par\noindent\hangindent 30pt\hbox to30pt{[#1]\hfil}}

\def\Brill{1}
\def\Brio{2}
\def\Cayley{3}
\def\Cayleytwo{4}
\def\Cayleythree{5}
\def\Des{6}
\def\Dodg{7}
\def\DW{8}
\def\Jacobi{9}
\def\LLT{10}
\def\Lasc{11}
\def\Lec{12}
\def\LP{13}
\def\Mert{14}
\def\Muirtreat{15}
\def\Muir{16}
\def\Muirtwo{17}
\def\Pf{18}
\def\RR{19}
\def\Saal{20}
\def\Scheib{21}
\def\Schur{22}
\def\Stem{23}
\def\Tan{24}
\def\Tor{25}
\def\Velt{26}
\def\Wenz{27}
\def\Zaje{28}
\def\Zajp{29}

\centerline{\bf Overlapping \Pfaff/ians}

\bigskip
\centerline{Donald E. Knuth, Computer Science Department, Stanford University}

\bigskip
\centerline{\sc To Dominique Cyprien Foata on his 60th birthday}

\bigskip
{\narrower\smallskip\noindent
{\bf Abstract.} A combinatorial construction proves an identity for the product
of the \Pfaff/ian of a skew-symmetric matrix by the \Pfaff/ian of one of its
submatrices. Several applications of this identity are followed by a brief
history of \Pfaff/ians.
\smallskip}

\bigskip\noindent
{\bf 0. Definitions.}
Let $X$ be a possibly infinite index set. We consider quantities $f[xy]$
defined on ordered pairs of elements of~$X$, satisfying the law of skew
symmetry
$$f[xy]=-f[yx]\,,\qquad {\rm for} \quad x,y\in X\,.\eqno(0.0)$$
This notation is extended to $f[\alpha]$ for arbitrary words $\alpha=x_1\ldots
x_{2n}$ of even length over~$X$ by defining the {\it \Pfaff/ian\/}
$$f[x_1\ldots x_{2n}]=\sum s(x_1\ldots x_{2n}, y_1\ldots y_{2n})\,
f[y_1y_2]\ldots f[y_{2n-1}y_{2n}]\,,\eqno(0.1)$$
where the sum is over all $(2n-1)(2n-3)\ldots (1)$ ways to write $\{x_1,\ldots
x_{2n}\}$ as a union of pairs $\{y_1,y_2\}\cup \cdots \cup
\{y_{2n-1},y_{2n}\}$, and where $s(x_1\ldots x_{2n},y_1\ldots y_{2n})$ is the
sign of the permutation that takes $x_1\ldots x_{2n}$ into $y_1\ldots y_{2n}$.

The \Pfaff/ian is well defined, even though there are $2^nn!$ different
permutations $y_1\ldots y_{2n}$ that yield the same partition
$\{y_1,y_2\}\cup\ldots\cup \{y_{2n-1},y_{2n}\}$ into pairs. For if we
interchange $y_{2j-1}$ with~$y_{2j}$, we change the sign of both $s(x_1\ldots
x_{2n}, y_1\ldots y_{2n})$ and $f[y_1,y_2]\ldots f[y_{2n-1}y_n]$, by (0.0); if
we interchange $y_{2i-1}$ with $y_{2j-1}$ and $y_{2i}$ with~$y_{2j}$, both
factors stay the same. Thus, for example,
$$\eqalignno{f[wxyz]
&=f[wx]f[yz]-f[wy]f[xz]+f[wz]f[xy]\cr
&=f[wx]f[yz]+f[wy]f[zx]+f[wz]f[xy]\,.&(0.2)\cr}$$

A partition into pairs is commonly called a {\it perfect matching}. Therefore
it is convenient to abbreviate (0.1) in the form
$$f[\alpha]=\sum_{\mu\in M(\alpha)} s(\alpha,\mu)\,\Pi f[\mu]\eqno(0.3)$$     
where $M(\alpha)$ is the set of perfect matchings of $\alpha$ represented as
words $y_1\ldots y_{2n}$ in some canonical way, and $\Pi f[y_1\ldots
y_{2n}]=f[y_1y_2]\ldots f[y_{2n-1}y_{2n}]$.

Notice that we have
$$f[wxyz]=-f[xyzw]\,.\eqno(0.4)$$
In general, an odd permutation of $\alpha$ will reverse the sign of
$f[\alpha]$, because every term in (0.3) changes sign.

\Pfaff/ians can also be defined recursively, starting with the null word
$\epsilon$ and proceeding to words of greater length:
$$\eqalignno{f[\epsilon]&=1\,;\cr
f[x_1\ldots x_{2n}]&=\sum_{j=2}^{2n}f[x_1x_j]f[x_{j+1}\ldots x_{2n}x_2\ldots
x_{j-1}]\,,\qquad n>0\,.&(0.5)\cr}$$
This recurrence
[\Jacobi]
corresponds to a procedure that constructs all perfect matchings
by starting with $\{x_1,x_2\}\cup\cdots\cup\{x_{2n-1},x_{2n}\}$ and
making cyclic permutations of the indices in positions
$\{2,\ldots,2n\}$, $\{4,\ldots,2n\}$, $\ldots\,$; each of these permutations
is even.

It will be convenient in the sequel to extend the sign function~$s$ to 
$s(\alpha,\beta)$ for arbitrary words $\alpha,\beta\in X^{\ast}$. We define
$s(\alpha,\beta)=0$ if either $\alpha$ or~$\beta$ has a repeated letter, or if
$\beta$ contains a letter not in~$\alpha$. Otherwise $s(\alpha,\beta)$ is the
sign of the permutation that takes $\alpha$ into the word
$$\beta\,(\alpha\blash\beta)\,,$$
where $\alpha\blash\beta$ is the word that remains when the elements
of~$\beta$ are removed from~$\alpha$. Thus, for example,
$$s(\alpha\beta\gamma,\beta)=\left\{\vcenter{\halign{$#$\hfil\quad%
&#\hfil\cr
0\,,&if $\alpha\beta\gamma$ contains a repeated letter;\cr
\noalign{\smallskip}
(-1)^{\vert\alpha\vert\,\vert\beta\vert}\,,&otherwise.\cr}}\right.\eqno(0.6)$$
We also have
$$s(\alpha,\beta\gamma)=s(\alpha,\beta)s(\alpha\blash\beta,\gamma)\,,
\eqno(0.7)$$
since both sides vanish unless the letters of $\beta\gamma$ are distinct and
contained in the distinct letters of~$\alpha$, and in the latter case
$s(\alpha,\beta\gamma)$ is the parity of the number of transpositions needed to
bring~$\beta$ to the left of~$\alpha$ and $\gamma$ to the left of the remaining
word~$\alpha\blash\beta$.

If $\alpha$ has repeated letters, the \Pfaff/ian $f[\alpha]$ is zero, because
$f[\alpha]=-f[\alpha]$ when we transpose two identical letters. Therefore our
convention that $s(\alpha,\beta)=0$ when $\alpha$ or $\beta$ has repeated
letters does not invalidate definition (0.1), which used a different convention
for $s(x_1\ldots x_{2n},y_1\ldots y_{2n})$. One consequence of the new
convention is the identity
$$f[\alpha]=\sum_{x_1<\cdots <x_n}\,\sum_{y_1>x_1}\cdots \sum_{y_n>x_n}
s(\alpha,x_1y_1\ldots x_ny_n)\,f[x_1y_1]\ldots f[x_ny_n]\eqno(0.8)$$
for any word $\alpha$ of length $2n$, assuming that $X$ is an ordered
set; the sum is over all conceivable
perfect matchings $\mu=x_1y_1\ldots x_ny_n$, but
$s(\alpha,\mu)$ is zero unless $\mu$ is a perfect matching of~$\alpha$.

\medskip\noindent
{\bf 1. The basic identity.}
The following identity due to H. W. L. Tanner
[\Tan]
can now be proved:
$$f[\alpha]\,f[\alpha\beta]=\sum_ys(\beta,xy)\,f[\alpha
xy]\,f[\alpha\beta\blash xy]\,,\qquad\hbox{for all $x\in\beta$}.\eqno(1.0)$$
This formula is vacuous when $\vert\beta\vert=0$ and trivial when
$\vert\beta\vert=2$, but when $\vert\beta\vert=4$
it says in particular that
$$\eqalignno{f[\alpha]\,f[\alpha wxyz]
&=f[\alpha wx]\,f[\alpha yz]-f[\alpha wy]\,f[\alpha xz]+f[\alpha wz]\,f[\alpha
xy]\cr
&=f[\alpha wx]\,f[\alpha yz]+f[\alpha wy]\,f[\alpha zx]+f[\alpha wz]\,f[\alpha
xy]\,.&(1.1)\cr}$$
We will demonstrate (1.0) by giving a combinatorial interpretation to each term
on the left and right sides of the equation, when the \Pfaff/ians are expanded
as sums over perfect matchings.

A typical term on the right of (1.0) is
$$s(\beta,xy)\,s(\alpha xy,\mu)\,s(\alpha\beta\blash xy,\nu)\,\Pi 
f[\mu]\,\Pi f[\nu]\,,\eqno(1.2)$$
where $x$ and $y$ are distinct elements of $\beta$, \
$\mu$ is a perfect matching
of $\alpha xy$, and $\nu$ is a perfect matching of $\alpha\beta\blash xy$.
Ignoring the sign for the moment, we can construct a graph by superimposing the
matchings~$\mu$ and~$\nu$. In this graph all vertices of~$\alpha$ have degree~2
because they are matched in both~$\mu$ and~$\nu$; all vertices of~$\beta$ have
degree~1.

There is a unique maximal path that starts at $y$ and uses edges from~$\mu$
and~$\nu$ alternately. This path ends at some element of~$\beta$, call it~$z$.
Let $\mu_1$ and $\nu_1$ be the edges of~$\mu$ and~$\nu$ on this path; let
$\mu_0$ and~$\nu_0$ be the other edges. Then we define corresponding matchings
$$\mu'=\mu_0\cup\nu_1\,,\qquad \nu'=\nu_0\cup\mu_1\,,\eqno(1.3)$$
which will be the key to establishing (1.0).

\medskip
{\it Case 1}, $z\neq x$. In this case $\vert\mu_1\vert=\vert\nu_1\vert$, since
the path from~$y$ starts with an element of~$\mu$ and ends with an element
of~$\nu$. Thus the matchings $\mu'$ and~$\nu'$ correspond to another term
on the right side of~(1.0); we will prove that this other term
cancels with~(1.2). Since $\mu''=\mu$ and $\nu''=\nu$, this
will set up an involution between cancelling terms.

We have
$$\Pi f[\mu]\,\Pi f[\nu]=\Pi f[\mu_0]\,\Pi f[\mu_1]\,\Pi f[\nu_0]\,\Pi 
f[\nu_1]=\Pi f[\mu']\,\Pi f[\nu']\,,\eqno(1.4)$$
so (1.2) will cancel with its counterpart if the signs differ. The sign of
(1.2)~is
$$s(\alpha xyz,\mu_0\mu_1 z)\,s(\alpha\beta,xy\nu_0\nu_1)\,,\eqno(1.5)$$
because $s(\beta,xy)=s(\alpha\beta,xy)$ and
$s(\alpha\beta,xy)\,s(\alpha\beta\blash xy,\nu)=s(\alpha\beta,xy\nu)$
by~(0.7). The sign of the permutation that takes $\mu_1z$ into $\nu_1y$ is the
same as the sign of the permutation that takes $y\nu_0\nu_1$ into
$z\nu_0\mu_1$, hence (1.5) equals
$$s(\alpha xyz,\mu_0\nu_1y)\,s(\alpha\beta,xz\nu_0\mu_1)\,.$$
But this is the negative of $s(\alpha
xzy,\mu_0\nu_1y)\,s(\alpha\beta,xz\nu_0\mu_1)$, the sign of the term that
corresponds to~$\mu'$ and~$\nu'$.

\medskip
{\it Case 2}, $z=x$. In this case we have $\vert\mu_1\vert =\vert\nu_1\vert+2$,
since $\mu_1$ includes both~$x$ and~$y$ while $\nu_1$ is contained in~$\alpha$.
It follows that $\mu'$ and~$\nu'$ are perfect matchings of~$\alpha$
and~$\alpha\beta$, respectively, so they define a typical term
$$s(\alpha,\mu')\,s(\alpha\beta,\nu')\,\Pi f[\mu']\,\Pi f[\nu']\eqno(1.6)$$
from the left side of (1.0). Conversely, every such term corresponds to
matchings~$\mu$ and~$\nu$ for a uniquely defined term (1.2) on the right. The
sign of this term,
$$s(\alpha xy,\mu_0\mu_1)\,s(\alpha\beta,xy\nu_0\nu_1)\,,$$
agrees with $s(\alpha,\mu')\,s(\alpha\beta,\nu')=s(\alpha xy,\mu_0\nu_1 xy)\,
s(\alpha\beta,\nu_0\mu_1)$, because the permutation that takes~$\mu_1$ into
$\nu_1xy$ has the same sign as the permutation that takes $xy\nu_0\nu_1$
into~$\nu_0\mu_1$.

\medskip\noindent
{\bf 2. Basic applications.}
The special case $\alpha=\epsilon$ of (1.0) reads
$$f[\beta]=\sum_ys(\beta,xy)\,f[xy]\,f[\beta\blash xy]\,,\qquad
\hbox{for all $x\in\beta$}.\eqno(2.0)$$
This is a mild generalization of the recurrence (0.5); it tells us how to
expand $f[\beta]$ with respect to any element of~$\beta$. We can get rid of the
constraint $x\in\beta$ by summing over all~$x$:
$$f[\beta]={1\over\vert\beta\vert}\,\sum_x\,\sum_y
s(\beta,xy)\,f[xy]\,f[\beta\blash xy]\,.\eqno(2.1)$$
Applying this rule to $f[\beta\blash xy]$ and repeating until words of
length~2 are reached yields a $\vert\beta\vert$-fold sum,
$$f[\beta]={1\over (2n)(2n-2)\ldots 2}\,\sum_{x_1}\cdots\sum_{x_{2n}}\,
s(\beta,x_1\ldots x_{2n})\,f[x_1x_2]\ldots f[x_{2n-1}x_{2n}]\,,\eqno(2.2)$$
when $\vert\beta\vert =2n$; this is, of course, the same as (0.8) when we
collect equal terms.

Now let $\alpha$ be a fixed word such that $f[\alpha]\neq 0$, and consider the
function
$$g(\beta)=f[\alpha\beta]/f[\alpha]\eqno(2.3)$$
on the words of $X$. Tanner's identity (1.0) tells us that
$$g(\beta)=\sum_y\,s(\beta,xy)\,g(xy)\,g(\beta\blash xy)\,,\qquad
\hbox{for all $x\in\beta$}.\eqno(2.4)$$
But this is the same relation as (2.0); so $g$ satisfies the \Pfaff/ian
recurrence
(0.5). Therefore any identity for \Pfaff/ians leads {\it a fortiori\/} to an
identity for~$g$. In particular, (0.3)~tells us that
$$g(\beta)=\sum_{\mu\in M(\beta)}s(\beta,\mu)\,\Pi g(\mu)\,,$$
which is equivalent to
$$f[\alpha]^{n-1}\,f[\alpha\beta]=\sum_{M(\beta)}s(\beta,x_1y_1\ldots x_ny_n)\,
f[\alpha x_1y_1]\ldots f[\alpha x_ny_n]\eqno(2.5)$$
when $\vert\beta\vert =2n$, where the sum is over all perfect matchings
$x_1y_1\ldots x_ny_n$ of~$\beta$. The special case $n=2$ appears in~(1.1).

We can also construct a dual formula by starting with a fixed $\alpha\beta$
such that $f[\alpha\beta]\neq 0$ and defining
$$h(\gamma)=s(\alpha\beta,\gamma)\,
f[\alpha\beta\blash\gamma]/f[\alpha\beta]\eqno(2.6)$$
on the words $\gamma$ contained in $\alpha\beta$. Then (1.0) yields
$$h(\beta)=\sum_ys(\beta,xy)\,h(\beta\blash xy)\,h(xy)\,,\qquad
\hbox{for all $x\in\beta$};\eqno(2.7)$$
so we can derive a companion to (2.5) in a similar fashion:
$$f[\alpha]\,f[\alpha\beta]^{n-1}=\sum_{M(\beta)}s(\beta,x_1y_1\ldots
x_ny_n)\,f[\alpha\beta\blash x_1y_1]\ldots f[\alpha\beta\blash x_ny_n]\,.
\eqno(2.8)$$
Identities (2.4) and (2.7) are the \Pfaff/ian analogs of theorems about
determinants that Muir called the Law of Extensible Minors and the Law
of Complementaries. (See [\Muirtreat], \S179 and \S98 in the original edition;
\S187 and \S179 in Metzler's revision.)

\medskip\noindent
{\bf 3. Applications to determinants.}
Determinants are the special case of \Pfaff/ians in which the index set is
bipartite with respect to~$f$, in the sense that $f[xy]=0$ when $x$ and~$y$
belong to the same part. It is convenient to imagine that the set of indices
consists of two disjoint parts~$X$ and~$\bar{X}$, so that $x$ belongs to~$X$ if
and only if $\bar{x}$ belongs to~$\bar{X}$, and $f[xy]=f[\bar{x}\,\bar{y}]=0$
for all $x,y\in X$. The independent quantities are now
$f[x\bar{y}]=-f[\bar{y}x]$; we can regard $X$ as a set of ``rows'' and
$\bar{X}$ as a set of ``columns,'' so that $f[x\bar{y}]$ is essentially an
element of the matrix~$f$. We use $f[x,y]$ as an alternative notation for
$f[x\bar{y}]$. In fact, when $\alpha$ and~$\beta$ are arbitrary words of~$X$ we
write
$$f[\alpha,\beta]=f[\alpha\bar{\beta}^R]\eqno(3.0)$$ for the determinant formed
from rows~$\alpha$ and columns~$\beta$. Here $\bar{\beta}^R$ stands for the
reverse complement of~$\beta$:
$$\overline{y_1y_2\ldots y_n}{}^R=\bar{y}_n\ldots \bar{y}_2\bar{y}_1\,.
\eqno(3.1)$$

Definition (3.0) agrees with the usual definition of determinants, when
$\vert\alpha\vert=\vert\beta\vert=n$, since the perfect matchings of
$\alpha\bar{\beta}^R$ that do not have vanishing products correspond to the
products
$$f[x_1\bar{y}_1]\ldots f[x_n\bar{y}_n]=f[x_1,y_1]\ldots
f[x_n,y_n]\,,\eqno(3.2)$$
where $\alpha=x_1\ldots x_n$ and $y_1\ldots y_n$ is a permutation of~$\beta$;
the corresponding sign $s(\alpha\bar{\beta}^R,x_1\bar{y}_1\ldots x_n\bar{y}_n)$
is just $s(\beta,y_1\ldots y_n)$, because the permutation that takes $x_1\ldots
x_n\bar{y}_n\ldots\bar{y}_1$ to $x_1\bar{y}_1\ldots x_n\bar{y}_n$ is even. For
example, we have
$$\eqalign{f[wx,yz]&=f[wx\bar{z}\bar{y}]\cr
&=f[wx]\,f[\bar{z}\bar{y}]-f[w\bar{z}]\,f[x\bar{y}]
+f[w\bar{y}]\,f[x\bar{z}]\cr
&=0-f[w,z]\,f[x,y]+f[w,y]\,f[x,z]\,,\cr}$$
the usual $2\times 2$ determinant
$$\openup1\jot \left\vert\matrix{%
f[w,y]&f[w,z]\cr
f[x,y]&f[x,z]\cr}\right\vert\,.$$

Theorem (1.0) immediately yields a corresponding identity for determinants,
when we apply these definitions:
$$f[\alpha,\beta]\,f[\alpha\gamma,\beta\delta]=\sum_y\,s(\gamma,x)\,s(\delta,y)
\,f[\alpha x,\beta y]\,f[\alpha\gamma\blash x,\beta\delta\blash y]\,,
\eqno(3.3)$$
for all $x\in\gamma$. When $\vert\gamma\vert=\vert\delta\vert$ is 2 or~3, this
identity reads
$$\eqalignno{f[\alpha,\beta]\,f[\alpha wx,\beta yz]
&=f[\alpha w,\beta y]\,f[\alpha x,\beta z]
 -f[\alpha w,\beta z]\,f[\alpha x,\beta y]\,;&(3.4)\cr
\noalign{\smallskip}
f[\alpha,\beta]\,f[\alpha uvw,\beta xyz]
&=f[\alpha u,\beta x]\,f[\alpha vw,\beta yz]\cr
&\quad\null -f[\alpha u,\beta y]\,f[\alpha vw,\beta xz]\cr
&\quad\null +f[\alpha u,\beta z]\,f[\alpha vw,\beta xy]\,.&(3.5)\cr}$$
Here are some small examples written in more conventional notation:
$$\openup-1\jot
\def\adet#1{\left\vert\,\vcenter{\halign{$a_{##}$&&\enspace$a_{##}$\cr
 #1}}\,\right\vert}
\eqalignno{a_{11}\adet{11&12&13\cr 21&22&23\cr 31&32&33\cr}
&=\adet{11&12\cr 21&22\cr}\;
 \adet{11&13\cr 31&33\cr}
-\adet{11&13\cr 21&23\cr}\;
 \adet{11&12\cr 31&32\cr}\,;&(3.6)\cr
\noalign{\bigskip}
\adet{11&12\cr 21&22\cr}\;
 \adet{11&12&13&14\cr 21&22&23&24\cr 31&32&33&34\cr 41&42&43&44\cr}
&=
 \adet{11&12&13\cr 21&22&23\cr 31&32&33\cr} \;
 \adet{11&12&14\cr 21&22&24\cr 41&42&44\cr}
-\adet{11&12&14\cr 21&22&24\cr 31&32&34\cr} \;
 \adet{11&12&13\cr 21&22&23\cr 41&42&43\cr}\,;\cr
\noalign{\vskip-2pt}
&&(3.7)\cr
\noalign{\medskip}
a_{11}
\adet{11&12&13&14\cr 21&22&23&24\cr 31&32&33&34\cr 41&42&43&44\cr}
&=
\adet{11&12\cr 21&22\cr}\;
\adet{11&13&14\cr 31&33&34\cr 41&43&44\cr}
-
\adet{11&13\cr 21&23\cr}\;
\adet{11&12&14\cr 31&32&34\cr 41&42&44\cr} \cr
\noalign{\smallskip}
&\hskip10em\null+
\adet{11&14\cr 21&24\cr}\;
\adet{11&12&13\cr 31&32&33\cr 41&42&43\cr}\,.&(3.8)\cr}$$

Of course determinants have been investigated rather thoroughly for nearly 250
years, so it would be surprising indeed if these identities were new.
Equation~(3.6) was, for instance, noted by Lagrange in 1773 
[\Muir, page 39];
(3.7) and higher examples of (3.4) were discussed by Desnanot in 1819
[\Muir, page 142].

One particularly interesting case in which (3.4) played a crucial role is C.~L.
Dodgson's elegant ``condensation method'' for determinant evaluation
[\Dodg],
discovered between the times when he wrote {\sl Alice in Wonderland\/} and {\sl
Through the Looking Glass\/}: Suppose the index set~$X$ is the integers, and
let $f_0[x,y]=1$ for all $x$ and~$y$, while $f_1[x,y]$ is the entry in row~$x$
and column~$y$ of a given matrix. Then for $k\geq 1$ let
$$f_{k+1}[x,y]=\left\vert\matrix{f_k[x,y]&f_k[x,y+1]\cr
\noalign{\smallskip}
f_k[x+1,y]&f_k[x+1,y+1]\cr}\right\vert /\,f_{k-1}[x+1,y+1]\,.\eqno(3.9)$$
It follows that
$$f_k[x,y]=f_1[x(x+1)\ldots (x+k-1), y(y+1)\ldots (y+k-1)]\qquad\hbox{for }
k\geq 0\,,\eqno(3.10)$$
by induction on $k$ using (3.4). To evaluate the $n\times n$ determinant
$f[1\,2\,\ldots\,n,\,1\,2\,\ldots\,n]$, we may therefore simply compute
$f_k[x,y]$ for $1\leq x,y\leq n+1-k$ and $k=2,\ldots,n$, hoping that it will
not be necessary to divide by zero. Dodgson's condensation method provided the
original motivation for Robbins and Rumsey's recent work on alternating sign
matrices
[\RR].

The earliest known identity involving products of determinants is
$$f[ab,12]\,f[ab,34]-f[ab,13]\,f[ab,24]+f[ab,14]\,f[ab,23]=0\,,\eqno(3.11)$$
which Alexis Fontaine des Bertins proudly wrote out 126 times for different
choices of the indices and then said ``et cetera.'' He submitted this and other
memoirs to the French academy in 1748, but the works remained unpublished until
1764 [\Muir, pp.\ 10--11]. From
(1.0) we can now recognize that the right-hand side of (3.8) is actually a
\Pfaff/ian product
$$f[ab]\,f[a\,b\,\bar{1}\,\bar{2}\,\bar{3}\,\bar{4}\,]\,,$$
which is indeed zero in the bipartite case. Bezout, in 1779, gave the similar
formula
$$\eqalignno{f[abc,123]&\,f[abc,456]-f[abc,124]\,f[abc,356]\cr
&\null+f[abc,125]\,f[abc,346]-f[abc,126]\,f[abc,345]=0\,,&(3.12)\cr}$$
and said ``on voit qu'il y a une infinit\'e d'autres combinaisons \`a faire''
[\Muir, page 51];
the right-hand side in this case is
$$f[a\,b\,c\,\bar{1}\,\bar{2}\,]\,f[a\,b\,c\,\bar{1}\,\bar{2}\,\bar{3}\,
\bar{4}\,\bar{5}\,\bar{6}\,]$$
when we replace determinants by \Pfaff/ians.

Another instance of (1.0) yields
$$\eqalignno{f[ab]\,f[a\,b\,c\,\bar{1}\,\bar{2}\,\bar{3}\,\bar{4}\,\bar{5}\,]
&=f[a\,b\,\bar{1}\,\bar{2}\,]\,f[a\,b\,c\,\bar{3}\,\bar{4}\,\bar{5}\,]
-f[a\,b\,\bar{1}\,\bar{3}\,]\;
f[a\,b\,c\,\bar{2}\,\bar{4}\,\bar{5}\,]\cr
&\quad\null+f[a\,b\,\bar{1}\,\bar{4}\,]\,
f[a\,b\,c\,\bar{2}\,\bar{3}\,\bar{5}\,]
-f[a\,b\,\bar{1}\,\bar{5}\,]\,f[a\,b\,c\,\bar{2}\,\bar{3}\,\bar{4}\,]\cr
&\qquad\null-f[a\,b\,\bar{1}\,c]\,f[a\,b\,\bar{2}\,\bar{3}\,\bar{4}
\,\bar{5}\,]\,.&(3.13)\cr}$$
Under bipartite restrictions this becomes an identity in determinants,
$$f[ab,12]\,f[abc,345]-f[ab,13]\,f[abc,245]+f[ab,14]\,f[abc,235]
-f[ab,15]\,f[abc,234]=0\,,\eqno(3.14)$$
which Desnanot 
[\Des]
seems to have known only in the special case
$$f[ab,12]\,f[abc,134]-f[ab,13]\,f[abc,124]+f[ab,14]\,f[abc,123]=0\eqno(3.15)$$
where column 1 $=$ column 5, although he knew the general result (3.3) 
[\Muir, page 145].

Thus we see that the single \Pfaff/ian identity (1.0) unifies a variety of
different-appearing determinant identities that arise when the indices are
given bipartite structure in different ways.

When identity (2.8) is specialized to determinants, it gives a formula for
minors of the adjugate of a matrix (i.e., determinants of cofactors):
$$\eqalignno{%
&f[\alpha,\beta]\,f[\alpha  x_1\ldots x_n,\beta u_1\ldots y_n]^{n-1}\cr
\noalign{\smallskip}
&\quad =
\left\vert\matrix{f[\alpha x_2\ldots x_n,\beta y_2\ldots y_n]&\ldots%
&f[\alpha x_2\ldots x_n,\beta y_1\ldots y_{n-1}]\cr
\vdots&&\vdots\cr
f[\alpha x_1\ldots x_{n-1},\beta y_2\ldots y_n]&\ldots
&f[\alpha x_1\ldots x_{n-1},\beta y_1\ldots y_{n-1}]\cr}\right\vert
   \,.&(3.16)\cr}$$
This general formula was first published by Jacobi in 1834, although special
cases had been found by Lagrange in 1773 and Minding in 1829
[\Muir, pp.\ 39, 197, 208--209].
The formula that corresponds to (2.5),
$$f[\alpha,\beta]^{n-1}\,f[\alpha x_1,\ldots x_n,\beta y_1\ldots y_n]
=\left\vert\matrix{f[\alpha x_1,\beta y_1]&\ldots&f[\alpha x_1,\beta y_n]\cr
\vdots&&\vdots\cr
f[\alpha x_n,\beta y_1]&\ldots&f[\alpha x_n,\beta y_n]\cr}\right\vert
\,,\eqno(3.17)$$
is simpler but was not discovered until Sylvester introduced a new viewpoint in
1851 
[\Muirtwo, pp.\ 60--61].

\bigskip\noindent
{\bf 4. Applications to closed forms.}
Let $g$ be the skew-symmetric Blaschke operator
$$g[xy]={x-y\over 1-xy}\,.\eqno(4.0)$$
Laksov, Lascoux, and Thorup [\LLT, (A.12.3)] and
John R. Stembridge
[\Stem, Proposition 2.3(e)]
independently discovered the remarkable identity
$$g[x_1x_2\ldots x_n]=\prod_{1\leq i<j\leq n}\;{x_i-x_j\over 1-x_ix_j}\,,
\qquad \hbox{$n$ even,}\eqno(4.1)$$
for which they gave ingenious but rather special-purpose proofs.

We can, however, prove (4.1) as a special case of more general theorem that
follows from a special case of (1.0):

\proclaim
Theorem. The identity
$$f[x_1\ldots x_n]=\prod_{1\leq i<j\leq n}\,f[x_ix_j]\eqno(4.2)$$
holds for all even\/ $n$ if and only if it holds for\/ $n=4$. 

\smallskip\noindent
{\bf Proof.}
If $n>4$ and the identity holds for smaller even
values of $n$, let $\alpha$ be any word of length $n-4$. Then
$$\eqalign{f[\alpha]\,f[\alpha wxyz]
&=f[\alpha wx]\,f[\alpha yz]-f[\alpha wy]\,f[\alpha xz]
+f[\alpha wz]\,f[\alpha xy]\cr
&=R(f[wx]\,f[yz]-f[wy]\,f[xz]+f[wz]\,f[xy])\cr
&=R\,f[wxyz]\cr
&=R\,f[wx]\,f[wy]\,f[wz]\,f[xy]\,f[xz]\,f[yz]\,,\cr}$$
where if $\alpha=x_1\ldots x_{n-4}$ the common factor $R$ is
$$\biggl(\,\prod_{1\le i<j\le n-4}f[x_ix_j]^2\biggr)\,
\biggl(\,\prod_{1\le i\le n-4}f[x_iw]\,f[x_ix]\,f[x_iy]\,f[x_iz]\biggr)\,.$$
Therefore
$$f[x_1\ldots x_{n-4}]\,f[x_1\ldots x_n]=
  f[x_1\ldots x_{n-4}]\prod_{1\le i<j\le n}f[x_ix_j]\,.$$
Equation (4.2) follows unless $f[x_1\ldots x_{n-4}]=0$.

If $f[y_1\ldots y_{n-4}]=0$ for all subwords $y_1\ldots y_{n-4}$ of $x_1\ldots
x_n$, then $f[x_1\ldots x_n]=0$ and again (4.2) holds.
Finally, if $y_1\ldots y_{n-4}$ is a
subword such that $f[y_1\ldots y_{n-4}]\neq 0$, there is a permutation
$y_1\ldots y_n$ of $x_1\ldots x_n$ for which our argument proves $f[y_1\ldots
y_n]=\prod_{1\leq i<j\leq n}f[y_iy_j]$. This establishes (4.2), because
permutations of the indices change the signs of both sides in the
same manner. \ \pfbox

\medskip
The theorem is of interest because it applies not only to (4.0)
but also to the simpler function
$$f[x_ix_j]={x_i-x_j\over c+x_i+x_j}\eqno(4.3)$$
when $c$ is any complex constant. Thus we obtain a more-or-less ``closed form''
(4.2) for the \Pfaff/ian of a new kind of matrix. (The special case $c=0$ was
previously noted by Schur [\Schur, \S36].)

In fact, the general function
$$f[x_ix_j]={x_i-x_j\over c+b(x_i+x_j)+a\,x_ix_j}\,,\qquad
b^2=ac\pm 1\,,\eqno(4.4)$$
also satisfies the necessary conditions; this expression includes both (4.0)
and~(4.3).

Are there other skew-symmetric rational functions of two variables that satisfy
$$\eqalignno{f[wx]\,f[yz]+&f[wy]\,f[zx]+f[wz]\,f[xy]\cr
&=f[wx]\,f[wy]\,f[wz]\,f[xy]\,f[xz]\,f[yz]\;\hbox{?}&(4.5)\cr}$$
One can, of course, replace $f[xy]$ by $f[r(x)r(y)]$ for any rational
function~$r$, so any solution of (4.5) implies an infinite class of
equivalent solutions. Alain Lascoux [\Lasc] has recently found strong
reasons for believing that there are no other solutions, up to changes
of variables.

When $f[xy]$ is a polynomial, an amusing closed form of a similar type was
noticed by G.~Torelli
[\Tor]:
Let $f_k[xy]=(x-y)^k$; then
$$f_{n-1}[x_1\ldots x_n]=(-1)^{n/2\choose 2}\,
\biggl(\prod_{k=0}^{n/2-1}{n-1\choose k}\biggr)\,
\prod_{1\leq i<j\leq n}(x_i-x_j)\eqno(4.6)$$
when $n$ is even. It is easy to prove this identity, as well as the fact that
$f_{2m-1} [x_1\ldots x_n]=0$ for $2m<n$, by observing that the \Pfaff/ian must
vanish when $x_i=x_j$.

\bigskip\noindent
{\bf 5. Generalization of the basic identity.}
Equation (1.0), which gives an expression for $f[\alpha]\,f[\alpha\beta]$ when
$\alpha$ is a proper subword of~$\alpha\beta$,
leads to a similar identity that is
useful when two words have an odd number of letters in common. Suppose
$\alpha\beta\gamma$ has no repeated letters, and let $x\in\beta$. Then
$$\eqalignno{f[\alpha\beta]\,f[\alpha\gamma]&=\sum_ys(\beta,xy)\,
f[\alpha\beta\blash xy]\,f[\alpha\gamma xy]\cr
&\qquad\null+\sum_y s(\beta,x)\,s(\gamma,y)\,f[\alpha y\beta\blash x]\,
f[\alpha x\gamma\blash y]\,.&(5.0)\cr}$$
For example, when $\vert\alpha\vert$ is odd we have
$$\eqalignno{f[\alpha xyz]\,f[\alpha uvw]
&=f[\alpha z]\,f[\alpha uvwxy]-f[\alpha y]\,f[\alpha uvwxz]\cr
&\qquad\null+f[\alpha uyz]\,f[\alpha xvw]-f[\alpha vyz]\,f[\alpha xuw]
+f[\alpha wyz]\,f[\alpha xuv]\,.&(5.1)\cr}$$
To prove (5.0), let $\gamma =x_1\ldots x_k$. We will construct a ``cancelling''
word $\gamma'=x'_k\ldots x'_1$ on new indices, by defining
$$f[yx'_j]=0\quad{\rm if}\quad y\neq x_j\,;\qquad f[x_jx'_j]=1\,.\eqno(5.2)$$
Then $f[\alpha\beta]=f[\alpha\gamma\gamma'\beta]$, and we can use (1.0) to
conclude that
$$f[\alpha\beta]\,f[\alpha\gamma]=
\sum_ys(\gamma'\beta,xy)\,f[\alpha\gamma\gamma'\beta\blash
xy]\,f[\alpha\gamma xy]\,.\eqno(5.3)$$
Now if $y\in\beta$ we have $s(\gamma'\beta,xy)=s(\beta,xy)$, and
$f[\alpha\gamma\gamma'\beta\blash xy]=f[\alpha\beta\blash xy]$. But if
$y=x'_j$ we have $s(\gamma'\beta,xy)=(-1)^js(\beta,x)$,
$f[\alpha\gamma\gamma'\beta\blash xy]=(-1)^{j-1}\,f[\alpha y\beta\blash
x]$, $f[\alpha\gamma xy]=(-1)^jf[\alpha x\gamma\blash y]$, and
$s(\gamma,y)=(-1)^{j-1}$.

\bigskip\noindent
{\bf 6. A brief history of \Pfaff/ians.}
Johann Friedrich \Pfaff/ introduced the functions that now bear his name in
1815 [\Pf]
[\Muir, pp.\ 396--401],
while studying a general way to solve systems of first-order partial
differential equations. He gave two procedures for listing all perfect
matchings, and observed that when the matchings are ordered lexicographically
the corresponding signs are strictly alternating $+,-,+,\ldots,+$. 

Jacobi developed \Pfaff/'s method further in 1827 
[\Jacobi],
and discovered an analog of ``Cramer's rule'' for the solution of general
systems of skew-symmetric linear equations
$$\sum_{j=1}^{2n}f[ij]\,z_j=f[i0]\,,\qquad n\hbox{ even}\,; \eqno(6.0)$$
namely,
$$z_j={f[1\ldots (j-1)\,0\,(j+1)\ldots n]\over f[1\ldots n]}\,.\eqno(6.1)$$
This implicitly proves that the \Pfaff/ian $f[1\ldots n]$ is a factor of the
general skew-symmetric determinant
$$\left\vert\matrix{f[11]&\ldots&f[1n]\cr
\vdots&&\vdots\cr
f[n1]&\ldots&f[nn]\cr}\right\vert\,,\qquad n\hbox{ even}.\eqno(6.2)$$
Cayley proved in 1849
[\Cayley]
that this determinant is in fact equal to the {\it square\/} of $f[1\ldots n]$.

An elegant graph-theoretic proof of Cayley's theorem, somewhat analogous to the
derivation of (1.0) above, was found by Veltmann in 1871
[\Velt]
and independently by Mertens in 1877
[\Mert].
Their proof anticipated 20th-century studies on the superposition of two
matchings, and the ideas have frequently been rediscovered. Cayley himself had
claimed that such a proof would be possible, after doing the calculations for
$n=4$ on the final page of a paper he wrote in 1861
[\Cayleythree].
But we should note that his original method was simpler. In fact, Cayley
originally
[\Cayley]
gave a short inductive proof of the more general formula
$$\advance\normalbaselineskip1pt
\left\vert\matrix{f[xy]&f[x2]&f[x3]&\ldots&f[xn]\cr
f[2y]&f[22]&f[23]&\ldots&f[2n]\cr
f[3y]&f[32]&f[33]&\ldots&f[3n]\cr
\vdots&\vdots&\vdots&&\vdots\cr
f[ny]&f[n2]&f[n3]&\ldots&f[nn]\cr}\right\vert
=f[x23\ldots n]\,f[y23\ldots n]\,,\eqno(6.3)$$
for arbitrary $x$ and $y$ when $n$ is even. 
And he proved several years later
[\Muirtwo, pp.\ 269, 278]
that the determinant on the left of (6.3) is $f[xy23\ldots n]\,f[23\ldots n]$
when $n$ is odd. (This determinant is incidentally {\it not\/} the same as
$f[x2\ldots n,y2\ldots n]$; the elements of the latter are
$f[x,y],f[x,2],\ldots
= f[x\bar{y}\,],f[x\bar{2}\,],\ldots\,$, not $f[xy],f[x2],\ldots\,$, according
to our conventions. Moreover, we generally use the notation $f[x,y]$ only when
we assume that $f[xy]=0$.)

It was Cayley who introduced the name {\it \Pfaff/ian}, because of its
``connexion with the researches of \Pfaff/ on differential equations''
[\Cayleytwo].
Another name {\it semideterminant\/}
(German {\it Halb\-determinant\/})  was proposed by Wilhelm Scheibner
[\Scheib],
but it did not gain many adherents.

Theorem (1.0) was discovered by Henry William Lloyd Tanner in 1878
[\Tan],
who gave inductive proofs for the cases $\vert\beta\vert=4$ and
$\vert\beta\vert=6$ from which proof schemata for higher cases could be
inferred. W{\l}adys{\l}aw Zajaczkowski found another proof shortly afterward
[\Zaje]
[\Zajp]
based  on Jacobi's determinant theorem (3.16).
The theorem was independently rediscovered in 1901 by J. Brill~[\Brill], who
found a still better proof. He first established the identity
$${n-1\choose k}f[x_1\ldots x_{2n}]=\sum_{1\le j_1<\cdots\,<j_{2k}\le2n}
s(x_1\ldots x_{2n},x_1\ldots x_{2k})f[x_1\ldots x_{2k}]f[x_1\ldots
x_{2n} \blash x_1\ldots x_{2k}]\eqno(6.4)$$
by induction on $k$; then he made the left side zero by setting $x_{2n}=x_1$.
A series of further steps led him to~(1.0). But the combinatorial proof in
section~1 above seems preferable to all three of these early approaches.

Identity (5.0) was recently discovered by Wenzel [\Wenz, Proposition~2.3],
and demonstrated via exterior algebra by Dress and Wenzel [\DW].

The fact that \Pfaff/ians are more fundamental than determinants, in the sense
that determinants are merely the bipartite special case of a general sum over
matchings, went unnoticed for a long time. The first person to observe that
every $n\times n$ determinant is a \Pfaff/ian was apparently Louis Saalsch\"utz
in 1908
[\Saal],
but the implicitly bipartite nature of his construction was not stated in
his paper; a modern reader sees it only with hindsight. Brioschi had found a
complicated way to express a $2n\times 2n$ determinant as a \Pfaff/ian, in 1856
[\Brio]:
If $A$ is any $2n\times 2n$ matrix and if $Q=
I_n\otimes\bigl({0\atop-1}\,{1\atop0}\bigr)$,
the determinant of~$A$ is the \Pfaff/ian of~$A^TQA$.

\Pfaff/ians continue to find numerous applications, for example in matching
theory 
[\LP]
and in the enumeration of plane partitions
[\Stem].
It should prove interesting to extend Leclerc's combinatorics of relations for
determinants [\Lec] to the analogous rules for \Pfaff/ians.

\medskip\noindent
{\bf Acknowledgements.}
Discussions with Lyle Ramshaw helped greatly to clarify my proof of (1.0).
Paul Algoet kindly corrected several typographical errors in my preprint.
Alain Lascoux referred me to [\LLT] and~[\Lec],
John Stembridge told me about [\Schur], and an anonymous referee
called my attention to [\DW]. I~also thank the editors for their patience.

\bigskip
\centerline{\bf References}
\bigskip

{\advance\baselineskip -1pt \parskip5pt plus .5pt minus .5pt

\def\bib[#1] {\par\noindent\hangindent\parindent\hbox to\parindent{[#1]\hfil}}

\bib
[\Brill]
J. Brill, ``Note on the algebraic properties of \Pfaff/ians,'' {\sl
Proceedings of the London Mathematical Society\/ \bf34} (1901),
143--151.

\bib
[\Brio]
F. Brioschi, ``Sur l'analogie entre une classe de d\'eterminants d'order pair;
et sur les d\'e\-ter\-mi\-nants  binaires,''
{\sl Journal f\"ur die reine und angewandte Mathematik\/
\bf 52} (1856), 133--141.

\bib
[\Cayley]
A. Cayley, ``Sur les d\'eterminants gauches,'' 
{\sl Journal f\"ur die reine und angewandte Mathematik\/
 \bf 38} (1849), 93--96. Reprinted in his {\sl Collected
Mathematical Papers\/ \bf 1}, 410--413.

\bib
[\Cayleytwo]
A. Cayley, ``On the theory of permutants,'' {\sl Cambridge and Dublin
Mathematical Journal\/ \bf 7} (1852), 40--51. Reprinted in his {\sl Collected
Mathematical Papers\/ \bf 2}, 16--26.

\bib
[\Cayleythree]
A. Cayley, ``Note on the theory of determinants,'' {\sl Philosophical
Magazine\/ \bf 21} (1861), 180--185. Reprinted in his {\sl Collected
Mathematical Papers\/ \bf 5}, 45--49.

\bib
[\Des]
P. Desnanot, {\sl Compl\'ement de la Th\'eorie des \'Equations du Premier
Degr\'e\/} (Paris, 1819).

\bib
[\Dodg]
C. L. Dodgson, ``Condensation of determinants, being a new and brief
method for computing their arithmetical values,'' {\sl Proceedings of the Royal
Society\/ \bf 84} (1866), 150--155. Reprinted in {\sl The Mathematical
Pamphlets of Charles Lutwidge Dodgson and Related Pieces}, edited by Francine
F.~Abeles (Charlottesville, Virginia: The University Press of Virginia, 1994),
170--180.

\bib
[\DW]
Andreas W. M. Dress and Walter Wenzel, ``A simple proof of an identity
concerning \Pfaff/ians of skew symmetric matrices,'' {\sl Advances in
Mathematics\/ \bf112} (1995), 120--134.

\bib
[\Jacobi]
C. G. Jacobi, ``\"Uber die {\it \Pfaff/}\kern.1emsche Methode,
eine gew\"ohnliche
line\"are Differential-gleichung zwischen $2n$ Variabeln durch ein System von
$n$ Gleichungen zu integriren,'' {\sl Journal f\"ur die reine und angewandte
Mathematik\/ \bf 2} (1827), 347--357. Reprinted in {\sl C.~G.~J. Jacobi's
Gesammelte Werke\/ \bf 4} (1886), 17--29.

\bib
[\LLT]
D. Laksov, A. Lascoux and A. Thorup, ``On Giambelli's theorem on complete
correlations,'' {\sl Acta Mathematica\/ \bf162} (1989), 143--199.

\bib
[\Lasc] Alain Lascoux, personal communication, 10 April 1995.

\bib
[\Lec] Bernard Leclerc, ``On identities satisfied by minors of a matrix,''
{\sl Advances in Mathematics\/ \bf100} (1993), 101--132.

\bib
[\LP]
L\'aszl\'o Lov\'asz and Michael D. Plummer, {\sl Matching Theory\/} (Budapest:
Akad\'emiai Kiad\'o, 1986); North-Holland Mathematics Studies {\bf 121}.

\bib
[\Mert]
F. Mertens, ``\"Uber die Determinanten, deren correspondirende Elemente
$a_{pq}$ und~$a_{qp}$ ent\-gegen\-gesetzt gleich sind,'' 
{\sl Journal f\"ur die reine und angewandte Mathematik\/
\bf 82} (1877), 207--211.

\bib
[\Muirtreat]
Thomas Muir, {\sl A Treatise on the Theory of Determinants\/} (London:
Macmillan, 1882). Revised and enlarged by William H. Metzler (London:
Longmans, Green, 1933; New York, Dover, 1960).

\bib
[\Muir]
Thomas Muir, {\sl The Theory of Determinants in the Historical Order of
Development\/} (London: MacMillan, 1906).

\bib
[\Muirtwo]
Thomas Muir, {\sl The Theory of Determinants in the Historical Order of
Development}, volume~2 (London: MacMillan, 1911).

\bib
[\Pf]
J. F. \Pfaff/, ``Methodus generalis, aequationes differentiarum partialium, nec
non aequationes differentiales vulgares, utrasque primi ordinis, inter
quotcunque variabiles, completi integrandi,'' {\sl Abhandlungen der
K\"oniglich-Preu{\ss}ischen Akademie der Wissenschaften zu Berlin},
Mathematische Klasse (1814--1815), 76--136.

\bib
[\RR]
David P. Robbins and Howard Rumsey, Jr., ``Determinants and alternating sign
matrices,'' {\sl Advances in Mathematics\/ \bf 62} (1986), 169--184.

\bib
[\Saal]
Louis Saalsch\"utz, ``Zur Determinanten-Lehre,''
{\sl Journal f\"ur die reine und angewandte Mathematik\/
\bf 134} (1908), 187--197.

\bib
[\Scheib]
W. Scheibner, ``\"Uber Halbdeterminanten,'' {\sl Berichte \"uber die
Verhandlungen der K\"oniglich S\"achsischen Gesellschaft der Wissenschaften zu
Leipzig\/ \bf 11} (1859), 151--159.

\bib
[\Schur]
J. Schur, ``\"Uber die Darstellung der symmetrischen und der alternierenden
Gruppe durch gebrochene lineare Substitutionen,'' {\sl Journal f\"ur die reine
und angewandte Mathematik\/ \bf139} (1911), 155--250. Reprinted in Issai
Schur, {\sl Gesammelte Abhandlungen\/ \bf1} (1973), 346--441.

\bib
[\Stem]
John R. Stembridge, ``Nonintersecting paths, \Pfaff/ians, and plane
 partitions,''
{\sl Advances in Mathematics\/ \bf 83} (1990), 96--131.

\bib
[\Tan]
H. W. Lloyd Tanner, ``A theorem relating to \Pfaff/ians,''
{\sl Messenger of Mathematics\/ \bf 8} (1878), 56--59.

\bib
[\Tor]
Gabriele Torelli, ``Quistione 64,'' {\sl Giornale di Matematiche\/ \bf 24}
(1886), 377.

\bib
[\Velt]
W. Veltmann, ``Beitr\"age zur Theorie der Determinanten,'' {\sl Zeitschrift
f\"ur Mathematik und Physik\/ \bf 16} (1871), 516--525.

\bib
[\Wenz]
Walter Wenzel, ``\Pfaff/ian forms and $\mit\Delta$-matroids,'' {\sl Discrete
Mathematics\/ \bf115} (1993), 253--266.

\bib
[\Zaje]
W. Zajaczkowski, ``A theorem relating to \Pfaff/ians,'' 
{\sl Messenger of Mathematics\/ \bf 10} (1880), 36--37.

\def\og#1{\leavevmode\vtop{\baselineskip0pt \lineskip0pt
  \lineskiplimit0pt \ialign{##\crcr\relax#1\cr
                          \hidewidth\kern.2em
           \dimen0=.0040ex \multiply\dimen0\fontdimen1\font
           \kern-.0156\dimen0`\hidewidth\cr}}}

\bib
[\Zajp]
W. Zaj\og{a}czkowski, ``O pewn\'ej w{\l}asno\'sci pfafianu,''
 {\sl Rozprawy i Sprawozdania z Posiedze\'n}, Wydzia{\l}u
 Matematyczno-Przyrodniczego Akademii Umiej\og{e}tno\'sci {\bf 7}
(Krakow, 1880), 67--74.

}

\end